\documentclass[reqno,12pt]{amsart}
\usepackage{amsmath} 
\usepackage{amsthm, amsfonts, amsgen, stmaryrd,enumerate}
\usepackage[english]{babel}
\usepackage{fullpage}
\usepackage[centertags]{amsmath}
\usepackage{indentfirst}
\usepackage{fancyhdr}
\usepackage[dvips]{graphicx}
\usepackage{psfrag}
\numberwithin{equation}{section}
\usepackage[latin1]{inputenc}

\usepackage{color}

\title[Hodge Laplacians on quantum spheres]{Examples of Hodge Laplacians on quantum spheres}
\date{7 September  2011}
\author{Alessandro Zampini}
\address{ Mathematisches Institut der L.M.U., Theresienstra\ss e 39,  D-80333 M\"unchen, Germany.}
 \email{zampini@math.lmu.de}

\newtheorem{theo}{Theorem}[section]

\newtheorem{rema}[theo]{Remark}

\newcommand{\nn}{\nonumber}

\newcommand{\dd}{{\rm d}}
\newcommand{\ca}{\mathcal{A}}

\newcommand{\cl}{\mathcal{L}}

\newcommand{\cu}{\mathcal{U}}        
\newcommand{\SU}{\mathrm{SU}_q(2)}  
\newcommand{\ASU}{\ca(\mathrm{SU}_q(2))}  
\newcommand{\sq}{\mathrm{S}^2_{q}}  
\newcommand{\Asq}{\ca(\mathrm{S}^2_{q})}  
\newcommand{\su}{\cu_q(\mathfrak{su}(2))}  

\newcommand{\eps}{\varepsilon}      
\newcommand{\co}[2]{#1_{(#2)}}      
\newcommand{\hs}[2]{\left\langle #1,#2\right\rangle}  

\newcommand{\lt}{{\triangleright}}    
\newcommand{\rt}{{\triangleleft}}
\newcommand{\IC}{{\mathbb C}} 
\newcommand{\IR}{{\mathbb R}} 
\newcommand{\IN}{{\mathbb N}} 
\newcommand{\IZ}{{\mathbb Z}} 
\DeclareMathOperator{\id}{id}       
\DeclareMathOperator{\U}{U}       

\newcommand{\figureheight}{8cm}
\newcommand{\putfig}[2]{\begin{figure}[htp]
        \special{isoscale c:/itex/texfig/#1.wmf, \the\hsize \figureheight}
        \vspace{\figureheight}
        \caption{#2}\label{fig:#1}
        \end{figure}}
\newcommand{\pictureheight}{4cm}
\newcommand{\putpicture}[2]{\begin{figure}[htp]
        \special{isoscale c:/itex/texfig/#1.wmf, \the\hsize \pictureheight}
        \vspace{\pictureheight}
        \caption{#2}\label{fig:#1}
        \end{figure}}

\newcommand{\beqa}{\begin{eqnarray}}
\newcommand{\eeqa}{\end{eqnarray}}
\newcommand{\beq}{\begin{equation}}
\newcommand{\eeq}{\end{equation}}

\newcommand{\G}{\mathfrak{G}}
\begin{document}

\thispagestyle{empty}

\begin{abstract}
Using a non canonical braiding over the 3d left covariant calculus we present a family of Hodge operators on $\SU$ and on its homogeneous quantum space $\sq$. 

\end{abstract}


\maketitle
\tableofcontents

\section{Introduction}	

The geometrical setting of this paper is well known. We equip the manifold of the quantum group $\SU$  with the left covariant three dimensional calculus from \cite{wor87}, and its quantum homogeneous space given by the standard Podle\'s sphere $\sq$ with the induced left covariant two dimensional calculus \cite{po92}. A frame bundle approach allowed in \cite{maj} to describe a Riemannian structure with a metric (arising from purely algebraic relations), Hodge operator and a Laplacian on $\sq$, which was coupled to the $q$-monopole connection of this $\U(1)$ Hopf fibration in \cite{lrz}. The question we analyse in this paper is the following: is it possible to introduce a suitable family of Hodge operators on $\SU$ {\it and} a meaningful projection procedure giving  Hodge operators on the $\sq$ sphere, as happens in classical geometry? This question had been addressed in the same setting in \cite{ale10} (and on the same spaces equipped with different calculi in \cite{lz10}) mainly using the formalism from \cite{kmt}, which relates  Hodge operators on a left covariant differential calculus to  the properties of a scalar product on it.  In this paper we merge this concept with results in \cite{hec2001}: the 3d calculus on $\SU$ has a (non canonical) braiding, and the whole exterior algebra is recovered in terms of the corresponding antisymmetriser operators.  In analogy to \cite{ale11}, we set the question of defining -- using scalar products -- Hodge operators whose square is diagonal and such that its spectrum shows a degeneracy fitting with the degeneracy of the spectra of the  antisymmetrisers of the calculus: moreover, this is the notion  we relate to that of of symmetry of a tensor on $\SU$. In order to clarify the meaning of this approach, sec.~\S\ref{ove:ss} describes it at lenght in the classical setting, while sec.~\S\ref{eSU} presents how this idea works for $\SU$, eventually proposing a family of Hodge operators on $\sq$.

\section{A classical path}
\label{ove:ss}

Consider a $N$-dimensional connected Lie group $G$ given as the real form of a complex connected Lie group. Its group manifold is parallelizable: the space of 1-forms $\Omega^1(G)$ is a free bicovariant $N$-dimensional  $\ca(G)$-bimodule on the basis of left (right) invariant $\{\phi^{a}\}$ ($\{\eta^{a}\}$) 1-forms; the associated first order differential calculus is given by $(\dd, \Omega^1(G)$) with $\dd h=(L_{a}h)\omega^{a}=(R_{a}h)\eta^a$ in terms of the action of the dual  left (right) invariant derivations $L_{a}$ ($\{R_{a}\}$) on $h\in\,\ca(G)$. The standard flip given\footnote{We denote $\Omega^{\otimes k}=\Omega^{1}(G)\otimes_{\ca(G)}\cdots\otimes_{\ca(G)}\Omega^1(G)$ and drop the overall obvious  dependence on $G$.}  on a basis by $\tau:\omega^{a}\otimes\omega^b\,\mapsto\,\omega^b\otimes\omega^a$ is a braiding on $\Omega^{\otimes 2} $; the corresponding  antisymmetriser operators $A^{(k)}\,:\,\Omega^{\otimes k}\,\to\,\Omega^{\otimes k}$  ($k\,\in\,\IN$) give the exterior algebra $\Omega^{\wedge}=(\oplus_{k=1}^{N}\,\Omega^{k}, \wedge)$ as 
$\Omega^{\otimes k}\,\supset\,\Omega^k\,=\,
(\mathrm{Range}\,A^{(k)})\,\simeq\,\Omega^{\otimes k}/\ker\,A^{(k)}$, with 
$$
\omega^{a_{1}}\wedge\ldots\wedge\omega^{a_{k}}\;=\;A^{(k)}(\omega^{a_{1}}\otimes\ldots\otimes\omega^{a_{k}})\;=\;\sum_{\pi\in\,S_{k}}(-1)^{\pi}\omega^{\pi(a_{1})}\otimes\ldots\otimes\omega^{\pi(a_{k})}$$
 where  $S_{k}$ is the set of permutations of $k$ elements. 
The differential calculus $(\Omega^{\wedge}, \dd)$ is given by equipping the exterior algebra with the unique consistent graded $\dd:\Omega^{\wedge k}\,\to\,\Omega^{\wedge k+1}$ derivative operator satisfying $\dd^2=0$ and a graded Leibniz rule. Every $\Omega^{\wedge k}$ is a bicovariant free $\ca(G)$-bimodule with $\dim\,\Omega^{\wedge k}=N!/(k!(N-k)!)$ and $\Omega^{\wedge k}\,=\,\emptyset$ for $k>N$. The antisymmetrisers have a completely degenerate spectral decomposition,
\beq
\label{spclA}
A^{(k)}(\omega^{a_{1}}\wedge\ldots\wedge\omega^{a_{k}})=k!(\omega^{a_1}\wedge\ldots\wedge\omega^{a_{k}}).
\eeq
We consider a non degenerate tensor $g\,:\,\Omega^1\times\Omega^1\,\to\,\ca(G)$, whose components we use to set an  $\ca(G)$-bimodule contraction $g\,:\,\Omega^{k}\times\Omega^k\,\to\,\ca(G)$ given on a  basis by
\beq
\label{clacon}
g(\omega^{a_{1}}\wedge\ldots\wedge\omega^{a_{k}},\omega^{b_{1}}\wedge\ldots\wedge\omega^{b_{k}})\,=\,\sum_{\pi, \pi^{\prime}\,\in\,S_{k}}(-1)^{\pi+\pi^{\prime}}\,\Pi_{j=1,\ldots, k}\,g(\omega^{\pi(a_j)},\omega^{\pi^{\prime}(b_{j})}),
\eeq
 and then a sesquilinear scalar product 
 $$
 \hs{\phi}{\phi^{\prime}}_{G}\,=\,(1/k!)g(\phi^*,\phi^{\prime});
 $$
  if  $\mu\,=\,m\,\theta\,=\,\mu^*$ is the volume form (with $\theta$ the Haar measure on $G$, $m\,\in\,\IC$), the equation
 \beq
 \label{claTop}
 \int_{\mu}\phi^{*}\wedge T(\phi^{\prime})\,=\,\hs{\phi}{\phi^{\prime}}_{G}
 \eeq
uniquely defines a bijective $T\,:\,\Omega^{k}\,\to\,\Omega^{N-k}$ with $T(1)=\mu, \;T(\mu)=m$. The following equivalence holds ($\phi,\phi^{\prime}\,\in\,\Omega^1$)
\beq
\label{claeq}
g(\phi,\phi^{\prime})\,=\,g(\phi^{\prime},\phi)\qquad\Leftrightarrow\qquad
T^2(\phi)=(-1)^{N-1}\{T^2(1)\}\phi.
\eeq
The tensor $g$ is symmetric if and only if the square  of $T$ acts  on $\Omega^1$ as a constant depending on the volume\footnote{It is true that the factor can be  arbitrary, and that $g$ is symmetric if and only if $T^2(\phi^a)=\zeta\,\phi^a$ with $0\,\neq\,\zeta\,\in\,\IC$.  The choice in \eqref{claeq} will give the possibility of the usual overall normalization.};  given a symmetric $g$ one has also that $T^{2}(\phi)=(-1)^{k(N-k)}\{T^2(1)\}\phi$ with $\phi\,\in\,\Omega^k$. But such an operator $T$ is not (yet) an Hodge operator: it has to be real, and  the reality condition comes as  the equivalence
\beq
\label{reacla}
g(\phi,\phi^{\prime})^*\,=\,g(\phi^{\prime*},\phi^{*})\qquad\Leftrightarrow\qquad T(\phi^{*})\,=\,(T(\phi))^*.
\eeq
Fixing the reality condition as a compatibility of the action of $T$ with the hermitian conjugation on $\Omega^1$ proves to be indeed sufficient to have $[T,^*]=0$ on the whole exterior algebra $\Omega^{\wedge}$. Such a symmetric and real operator $T$ is then recovered as the Hodge operator corresponding to the (inverse) of the (metric) tensor $g$ on the group manifold. The choice  $\hs{\mu}{\mu}_{G}\,=\,{\rm sgn}(g)$ fixes the modulus of the scale parameter $m$ of the volume so to have 
$$
T^2(\phi)\,=\,{\rm sgn}(g)(-1)^{k(N-k)}\phi
$$ 
for any $\phi\,\in\,\Omega^k$. 

Let $K\,\subset\,G$ be a compact subgroup of $G$. The quotient of its right action gives a principal fibration $\pi\,:\,G\,\to\,G/K$. The exterior algebra $\Omega(G/K)\,\subset\Omega(G)$ is given in terms of horizontal and right $K$-invariant forms on $G$; any $\Omega^{s}(G/K)$ turns out to be no longer a free $\ca(G/K)$-bimodule. If we consider right $K$-invariant metric tensors $g$ on $G$ whose restriction to the homogeneus space $G/K$ is non degenerate, then  the restriction of the  corresponding scalar product consistently sets 
\beq
\label{hoho}
\hs{\psi}{\psi^{\prime}}_{G/K}\,=\,\hs{\psi}{\psi^{\prime}}_{G}\,=\,\int_{\check{\mu}}\psi^{*}\wedge\check{T}(\psi^{\prime})
\eeq
as a definition for the operator $\check{T}\,:\,\Omega^{s}(G/K)\,\to\,\Omega^{N^{\prime}-s}(G/K)$, with $\check{\mu}\,=\,\check{\mu}^*$ a volume form for the $N^{\prime}$-dimensional differential calculus on 
$G/K$. One has $$\check{T}^{2}(\psi)\,=\,(-1)^{s(N^{\prime}-s)}\hs{\check{\mu}}{\check{\mu}}_{G/K}\,\psi,$$ for any $\psi\,\in\,\Omega^{s}(G/K)$ and this induces to consider $\check{T}$, after a natural normalisation,  as the Hodge operator on $G/K$ obtained by projecting the  right $K$-invariant (inverse) metric  tensor $g$ onto the homogeneous space.

\section{Exterior algebras and Hodge operators over  quantum spheres}
\label{eSU}

We move to the quantum setting  with a short presentation of  $\SU$ and its well known 3d left covariant exterior algebra, referring the reader to \cite{wor87, ks} for a more complete description. 
As quantum group $\SU$ we consider the  polynomial  unital  $*$-algebra $\ASU=(\SU,\Delta,S,\eps)$ generated by elements $a$ and $c$ such that, using the matrix notation
\beq
\label{Us}
 u = 
\left(
\begin{array}{cc} a & -qc^* \\ c & a^*
\end{array}\right) , 
\eeq
the Hopf algebra structure can be expressed as
$uu^*=u^*u=1,\quad \,\Delta\, u = u \otimes u ,  \quad S(u) = u^* , \quad \eps(u) = 1.$
The deformation parameter
$q\in\IR$ is restricted  without loss of generality  to the interval $0<q<1$. 
The quantum universal envelopping algebra $\su$ is the unital Hopf $*$-algebra
generated the four elements $K^{\pm 1},E,F$, with $K K^{-1}=1$ and the
relations: 
\beq 
K^{\pm}E=q^{\pm}EK^{\pm}, \qquad 
K^{\pm}F=q^{\mp}FK^{\pm}, \qquad  
[E,F] =\frac{K^{2}-K^{-2}}{q-q^{-1}} . 
\label{relsu}
\eeq 
The $*$-structure is
$K^*=K, \,  E^*=F $,
while the Hopf algebra structures are 
\begin{align*}
&\Delta(K^{\pm}) =K^{\pm}\otimes K^{\pm}, \quad
\Delta(E) =E\otimes K+K^{-1}\otimes E,  \quad 
\Delta(F)
=F\otimes K+K^{-1}\otimes F,
\\ &\qquad S(K) =K^{-1}, \quad
S(E) =-qE, \quad 
S(F) =-q^{-1}F \\ 
&\qquad\qquad\varepsilon(K)=1, \quad \varepsilon(E)=\varepsilon(F)=0.
\end{align*}
The non degenerate Hopf algebra pairing between the two algebras above is  the $*$-compatible bilinear mapping $\hs{~}{~}:\su\times\ASU\to\IC$   given on the generators by 
\beq
\langle K^{\pm},a\rangle=q^{\mp 1/2}, \qquad  \langle K^{\pm},a^*\rangle=q^{\pm 1/2},  \qquad \langle E,c\rangle=1, \qquad \langle F,c^*\rangle=-q^{-1}, \label{ndp}
\eeq
with all other couples of generators pairing to zero. 
The $*$-compatible canonical commuting actions of $\su$ on $\ASU$ are (adopting the  Sweedler notation)
$$
h \lt x := \co{x}{1} \,\hs{h}{\co{x}{2}}, \qquad
x \rt  h := \hs{h}{\co{x}{1}}\, \co{x}{2}. 
$$
Given the algebra $\ca(\U(1)):=\IC[z,z^*] \big/ \!\!<zz^* -1>$,  the map  
\beq \label{qprp}
\pi: \ASU \, \to\,\ca(\U(1)),\quad\quad\pi(a)\,=\,z, \quad\pi(a^*)\,=\,z^*,\quad\pi(c)\,=\,\pi(c^*)\,=0
\eeq 
is a surjective Hopf $*$-algebra homomorphism, so that  $\U(1)$
is a quantum subgroup of $\SU$ with right coaction:
\beq 
\delta_{R}:= (\id\otimes\pi) \circ \Delta \, : \, \ASU \,\to\,\ASU \otimes
\ca(\U(1)) . \label{cancoa} 
\eeq 
This right coaction gives a decomposition 
\beq
\label{clnm}
\ASU=\oplus_{n\in\,\IZ}\cl_{n}, \qquad\qquad\cl_{n}:=\{x\in\,\ASU\;:\;\delta_{R}(x)=x\otimes z^{-n}\},
\eeq
with $\Asq=\cl_{0}$ the algebra of the standard Podle\'s sphere, each $\Asq$-bimodule $\cl_{n}$ giving the set of (charge $n$)  $\U(1)$-coequivariant maps for the topological quantum principal bundle $\Asq\hookrightarrow\ASU$.

The FODC $(\dd, \Gamma)$  we are going to consider on  $\SU$ was introduced by Woronowicz \cite{wor87}; $\Gamma$ is a three dimensional left covariant free $\ASU$-bimodule.
A basis of left-invariant 1-forms ($\Gamma_{{\rm inv}}\,\subset\,\Gamma$) is 
\beq
\omega_{z} =a^*\dd a+c^*\dd c , \qquad\qquad
 \omega_{-} =c^*\dd
a^*-qa^*\dd c^*, \qquad\qquad
\omega_{+} =a\dd c-qc \dd a.
\label{q3dom}
\eeq
It is  a $*$-calculus, with  $\omega_{-}^*=-\omega_{+}$ and
$\omega_{z}^*=-\omega_{z}$, while the bimodule structure is  $\omega_{z} x=q^{2n} x\,\omega_{z}, \quad
\omega_{\pm} x=q^{n} x\,\omega_{\pm}$ for any $x\,\in\,\cl_{n}$. 
The  basis of the quantum tangent space  $\mathcal{X}_{\mathcal{Q}}\subset\su$ dual to $\Gamma_{{\rm inv}}$ \eqref{q3dom}  
 is given by
\beq
\label{Xv3}
X_{-}=q^{-1/2}FK, \qquad X_{+}=q^{1/2} EK, \qquad X_{z}= (1-q^{-2})^{-1}\left(1-K^4\right).
\eeq
so that one can write exact 1-forms as $\dd x=\sum_{a=\pm,z}(X_{a}\lt x)\omega_{a}$
with $x\,\in\,\ASU$. 
This FODC has a non canonical braiding $\sigma$ on $\Gamma^{\otimes 2}$ \cite{hec2001}, which is as follows 
\begin{align}
\sigma(\omega_{a}\otimes\omega_{a})=\omega_{a}\otimes\omega_{a},&\qquad\qquad a=\pm,z
\nn \\
\sigma(\omega_{-}\otimes\omega_{+})=(1-q^{2})\omega_{-}\otimes\omega_{+}+q^{-2}\omega_{+}\otimes\omega_{-}, & \qquad\qquad \sigma(\omega_{+}\otimes\omega_{-})=q^4\omega_{-}\otimes\omega_{+}, \nn \\
\sigma(\omega_{-}\otimes\omega_{z})=(1-q^{2})\omega_{-}\otimes\omega_{z}+q^{-4}\omega_{z}\otimes\omega_{-}, &
\qquad\qquad \sigma(\omega_{z}\otimes\omega_{-})=q^6\omega_{-}\otimes\omega_{z}, \nn \\
\sigma(\omega_{z}\otimes\omega_{+})=(1-q^{2})\omega_{z}\otimes\omega_{+}+q^{-4}\omega_{+}\otimes\omega_{z},& \qquad\qquad
\sigma(\omega_{+}\otimes\omega_{z})=q^6\omega_{z}\otimes\omega_{+}.
\label{brai}
\end{align}
The corresponding\footnote{Explicitly, one has $A^{(2)}=1-\sigma$ acting on $\Gamma^{\otimes 2}$ and $A^{(3)}=(1-\sigma_{2})(1-\sigma_1+\sigma_1\sigma_2)$ on $\Gamma^{\otimes 3}$ where $\sigma_1=(\sigma\otimes 1)$ and $\sigma_2=(1\otimes \sigma)$.} antisymmetrisers $A^{(k)}:\Gamma^{\otimes k}\to\Gamma^{\otimes k}$ induce an exterior algebra $\Gamma^{\wedge}=\oplus_{k=1}^3\Gamma^{\wedge k}$ where $\Gamma^{\wedge k}={\rm Range}\,A^{(k)}$ and a  differential calculus $(\Gamma^{\wedge}, \dd)$, where each $\Gamma^{\wedge k}$ is a free left covariant $\ASU$-bimodule with the same classical dimension $\dim\Gamma^{\wedge k}=3!/(k!(3-k)!)$. One also has the isomorphism 
$(\Gamma^{\wedge}, \dd)\simeq(\Gamma^{\wedge}_{u},\dd)$, where the last is the universal differential (higher order) calculus over the FODC $(\Gamma,\dd)$ (see \cite[\S 12.2, \S14.3]{ks}).  In particular, the 
exterior wedge product satisfies
$\omega_{a}\wedge\omega_{a}=0, \quad (a=\pm,z)$ and  
\beq
\label{commc3}
\omega_{-}\wedge\omega_{+}+q^{-2}\omega_{+}\wedge\omega_{-}=0,
\qquad\qquad
\omega_{z}\wedge\omega_{\mp}+q^{\pm4}\omega_{-}\wedge\omega_{z}=0, 
\eeq
while 
the unique left invariant hermitian  top form of the calculus is $\theta=\theta^*=\omega_{-}\wedge\omega_{+}\wedge\omega_{z}$ with $x\,\theta=q^{-4n}\theta\,x$ for $x\in\cl_{n}$. Antisymmetrisers are totally degenerate, $A^{(k)}\omega\,=\,\lambda_{k}\omega$ with eigenvalues on $\Gamma^{\wedge 2}, \,\Gamma^{\wedge 3}$, 
\beq
\lambda_{2}\,=\,(1+q^2),\qquad\qquad
\lambda_{3}\,=\,(1+2q^2+2q^4+q^6).
\label{spA}
\eeq

\subsection{Symmetric contractions and Hodge operators}

We are going to mimick the classical construction  to obtain a suitable class of  scalar products in the quantum setting. We shall consider  $\ASU$-left invariant contractions, namely $\IC$-bilinear  maps $g\,:\,\Gamma_{{\rm inv}}\times\Gamma_{{\rm inv}}\,\to\,\IC$,   which are also right $\U(1)$-invariant. The right $\U(1)$-coaction \eqref{cancoa} consistently extends to a $\U(1)$-coaction $\delta_{R}^{(k)}$ on $\Gamma^{\wedge k}$, reading $\delta^{(k)}(\omega_{a})\,=\,\omega_{a}\otimes z^{n_a}$ on 1-forms\footnote{Notice that the restriction of the braiding \eqref{brai} to $\Gamma_{{\rm inv}}^{\otimes 2}$ commutes with this coaction, so the $\U(1)$ grading of higher order forms is given by the sum of their degree 1 components.}   with  $n_{\pm}=\pm2,\;
n_{z}=0$, so we set $g(\omega_{a},\omega_{b})=0$ if $n_{a}+n_{b}\neq0$. The only non zero coefficients are (non degeneracy clearly amounts to $\alpha\,\beta\,\gamma\,\neq\,0$)
\beq
\label{g1fo}
g(\omega_{-},\omega_{+})=\alpha, \qquad g(\omega_{+},\omega_{-})=\beta,\qquad g(\omega_{z},\omega_{z})=\gamma.
\eeq
This contraction is naturally extended to $\Gamma_{{\rm inv}}^{\wedge 2},\,\Gamma^{\wedge3}_{{\rm inv}}$ by the quantum analogue of the classical \eqref{clacon}, replacing the action of the classical antiysmmetrisers by the quantum ones $A^{(2)},\,A^{(3)}$. 
We set  a left-invariant scalar product on the whole exterior algebra $\Gamma^{\wedge}$ by using the spectra of the antisymmetrisers  to normalise the contractions and then introduce a  corresponding {\it dual} operator $T\,:\,\Gamma^{\wedge k}\,\to\,\Gamma^{\wedge 3-k}$,
\beq
\label{quT}
\hs{x\,\omega}{x^{\prime}\omega^{\prime}}_{\SU}\,=\,h(x^*x^{\prime})\,\frac{1}{\lambda_{k}}\,g(\omega^*,\omega^{\prime})\,=\,\int_{\mu}(x\,\omega)^*\wedge T(x^{\prime}\omega^{\prime})
\eeq
where $\omega,\,\omega^{\prime}\,\in\,\Gamma^{\wedge k}_{{\rm inv}}, \quad x,x^{\prime}\,\in\,\ASU$,  $h$ is the Haar state of $\SU$, $\mu=m\,\theta,\,m\,\in\,\IR$ is an hermitian volume form. The operator $T$ is well defined (provided $g$ is non degenerate) and left $\ASU$-linear \cite{kmt}.  From  the spectrum of its square we define $g$: 
\begin{enumerate}[i)]
\item {\it symmetric}, provided the operator $T^2\,:\,\Gamma\,\to\,\Gamma$ is constant (recall $\Gamma=\Gamma^{\wedge 1}$),
\item {\it real}, provided the relation $T(\omega_{a}^*)\,=\,(T(\omega_{a}))^*$ for $\omega_{a}\,\in\,\Gamma_{{\rm inv}}$ holds.
\end{enumerate}
These two conditions are compatible, the set $\G_{\sigma}$ of $\ASU$-left invariant and right $\U(1)$-invariant symmetric and real non degenerate  contractions $g$ is given by \eqref{g1fo} with $\beta\,=\,q^6\alpha\,\in\,\IR,\; \gamma\,\in\,\IR,\; \alpha\,\gamma\neq0$.  For $g\,\in\,\G_{\sigma}$ we assume $T$ to be a {\it quantum Hodge} operator with  $T(1)=\mu$, 
\begin{align}
T(\omega_{-})\,=\,-q^{-2}m\,\beta(\omega_{-}\wedge\omega_{z}), &\qquad\qquad T(\omega_{-}\wedge\omega_{z})\,=\,2q^{-4}(\lambda_2)^{-1}m\,\beta\,\gamma\,\omega_{-}, \nn \\
T(\omega_{+})\,=\,m\,\alpha(\omega_{+}\wedge\omega_{z}), &\qquad\qquad T(\omega_{+}\wedge\omega_{z})\,=\,-2q^{6}(\lambda_2)^{-1}m\,\alpha\,\gamma\,\omega_{+}, \nn \\ 
T(\omega_{z})\,=\,m\,\gamma(\omega_{-}\wedge\omega_{+}), &\qquad\qquad T(\omega_{-}\wedge\omega_{+})\,=\,-2(\lambda_2)^{-1}m\,\alpha\,\beta\,\omega_{z}
\label{T23}
\end{align}
and $T(\theta)\,=\,-6q^4(\lambda_{3})^{-1}m\,\alpha\,\beta\,\gamma$. We define the quantum determinant  of a contraction as  $\det\,g=\hs{\mu}{\mu}_{\SU}$, with signature  ${\rm sgn}(g)= {\rm sgn}(\det\,g)$,  and  fix the modulus of the scalar factor $m$ as $T^2(1)\,=\,{\rm sgn}(g)$. 
For $g\,\in\,\G_{\sigma}$ it is ${\rm sgn}(g)=-{\rm sgn}(\gamma)$:  this normalisation condition gives $T^2(\omega_{a})\,=\,{\rm sgn}(g)(2\lambda_{3})(6q^4\lambda_2)^{-1}\omega_{a}$. The spectra of the square of the quantum Hodge has the same degeneracy of the quantum antisymmetrisers of the given calculus, their specific eigenvalue appears as a deformation of the classical eigenvalue ${\rm sgn}(g)$ with $g$ the classical metric on ${\rm SU(2)}$.

For $g\,\in\,\G_{\sigma}$ the Laplacian operator on  $\SU$ is (see \eqref{Xv3}):
\beq
\Box_{\SU}x\,=\,\{\alpha(X_{-}X_{+}\,+\,q^6X_{+}X_{-})\,+\,\gamma\,X_{z}X_{z}\}\lt x
\label{l3s}
\eeq
with $x\,\in\,\ASU$. Its spectrum $\sigma_{P}$ is discrete (see \cite[eq. 4.23]{ale10}); one has  $\sigma_{P}\,\subset\,\IR^+$ if $\alpha\,\gamma\,>\,0$, while $\sigma_{P}$ is unbounded from above and below if $\alpha\,\gamma\,<\,0$. The formalism developed here gives then Laplacians whose classical limit recover both a Riemannian and a not Riemannian case.

\begin{rema}
As usual in the quantum setting, for the braiding \eqref{brai} on $\Gamma^{\otimes 2}$  it is  $\sigma^{2}\neq1$; the map $\sigma^{-1}$ gives antisymmetrisers $A^{(k)}_{\sigma^{-1}}$ and a differential calculus  $(\dd, \Gamma_{\sigma^{-1}}^{\wedge})$ on $\SU$  isomorphic to $(\dd, \Gamma^{\wedge})$.  One has $\G_{\sigma}\,=\,\G_{\sigma^{-1}}$.  The two families of Hodge operators have isomorphic spectra  (and give equivalent Laplacians as well), but are  not compatible (unlike the ones constructed in \cite{ale11}) since $[T,T_{\sigma^{-1}}]\neq0$ on $\Gamma^{\wedge}$.

\end{rema}

From \cite{maj,lrz} we know how to write the 2 dimensional exterior algebra $\Omega(\sq)$  induced on the Podle\'s sphere as $\Asq$-bimodules (frame bundle approach),
\beq
\label{fbsq}
\Omega({\sq}) = \Asq \oplus \left(\cl_{-2} \omega_{-}
\oplus \cl_{+2} \omega_{+} \right) \oplus \Asq \omega_{-}\wedge\omega_{+};
\eeq
from \cite[\S 5]{ale10} we know how the above scalar product restricts to an inner product on $\Omega(\sq)$, and how to compute the corresponding {\it dual} $\check{T}\,:\,\Omega^{s}(\sq)\to\Omega^{2-s}(\sq)$, which is left $\Asq$-linear, in complete analogy to the classical construction.  We obtain $\check{T}(1)=i\check{m}\,\omega_{-}\wedge\omega_{+}, \quad \check{T}(\omega_{-}\wedge\omega_{+})=2\check{m}^{2}(\lambda_2)^{-1}$ and
\beq
\check{T}(v_{-}\,\omega_{-})=-iq^{-2}\check{m}\beta\,v_{-}\,\omega_{-},\qquad\qquad\check{T}(v_{+}\,\omega_{+})=i\check{m}\alpha\,v_{+}\,\omega_{+}
\label{ctsq}
\eeq
with $v_{\pm}\,\in\,\cl_{\pm2}$. Upon the natural normalisation $\check{m}^2=\lambda_{2}\alpha\,\beta/2$, we see that the squared operator $\check{T^2}$ {\it is} diagonal, but {\it not} constant on $\Omega^1(\sq)$ for $g\,\in\,\G_{\sigma}$.  

Let us briefly comment on that. This should not be too surprising, but  even expected. Expression \eqref{fbsq} shows explicitly that $\Omega(\sq)$ is {\it not} a free left covariant $\Asq$-bimodule, and {\it cannot} be given as the range of a suitable family of antisymmetriser operators on $\Omega^{\otimes 2}(\sq)$. 
This means  -- recall that we link the consistency of an Hodge operator to the equivalence between the 
{\it degeneracy} of the spectrum of its square  and  that of the antisymmetrisers for the calculus --  that we actually have no natural term of comparison. 

The operator $\check{T}$  appears  nevertheless as the natural analogue of the Hodge operator on the classical sphere ${\rm S}^2\sim{\rm SU}(2)/\U(1)$, so we propose to assume it as a quantum Hodge  on $\sq$. If we define the corresponding  Laplacians  on $\sq$, then it is evident that we can {\it verbatim} repete the analysis in \cite[\S 6]{lrz} and \cite[\S 7]{ale10}, since the setting we have described in this section gives the geometry  \cite{bm} of 
 a quantum $\U(1)$ Hopf fibration over $\sq$ with compatible calculi. 

\section*{Acknowledgments}
I should like to thank Giovanni Landi for the impulse he gave to this research line, Istvan Heckenberger for the precious feedback,  the  organisers of the FUNINGEO conference 2011, Detlef D\"urr and the Mathematisches Institut at  L.M.U. for the support. 
This paper is dedicated to Beppe, my  {\it Doktorvater}:  so often I remember what he used to say  during his classes..."Instead of solving three different exercises, try  to solve the same, in three different ways". {\it Danke, herzlichst.}

\end{document}